\documentclass{amsart}
\usepackage{amssymb,amsmath,amscd,amsthm}
\usepackage{amsfonts,epsfig,latexsym,graphicx,amssymb}

\newtheorem{Thm}{Theorem}
\newtheorem{Lm}{Lemma}[section]
\newtheorem{Prop}{Proposition}
\newtheorem{Def}[Lm]{Definition}
\newtheorem{Rem}[Lm]{Remark}
\newtheorem{Cor}{Corollary}
\newtheorem*{GL}{Gap Lemma}

\def\R{\mathbb{R}}

\def\Z{\mathbb{Z}}

\def\T{\mathbb{T}}

\def\bdef{\begin{Def}}
\def\endef{\end{Def}}
\def\bthm{\begin{Thm}}
\def\ethm{\end{Thm}}
\def\bprop{\begin{Prop}}
\def\enprop{\end{Prop}}
\def\blm{\begin{Lm}}
\def\elm{\end{Lm}}
\def\bcor{\begin{Cor}}
\def\ecor{\end{Cor}}
\def\brm{\begin{Rem}}
\def\erm{\end{Rem}}
\def\bfig{\begin{picture}}
\def\efig{\end{picture}}
\def\beq{\begin{eqnarray}}
\def\eneq{\end{eqnarray}}
\def\beal{\begin{aligned}}
\def\enal{\end{aligned}}

\title{The Spectrum of the Weakly Coupled Fibonacci Hamiltonian}

\author{David Damanik}

\address{Department of Mathematics, Rice University, Houston, TX~77005, USA}

\email{damanik@rice.edu}

\author{Anton Gorodetski}

\address{Department of Mathematics, University of California, Irvine CA 92697, USA}

\email{asgor@math.uci.edu}

\thanks{D.\ D.\ was supported in part by NSF grant DMS--0800100.}

\date{\today}

\begin{document}

\begin{abstract}
We consider the spectrum of the Fibonacci Hamiltonian for small
values of the coupling constant. It is known that this set is a
Cantor set of zero Lebesgue measure. Here we study the limit, as the
value of the coupling constant approaches zero, of its thickness and
its Hausdorff dimension. We announce the following results and
explain some key ideas that go into their proofs. The thickness
tends to infinity and, consequently, the Hausdorff dimension of the
spectrum tends to one. Moreover, the length of every gap tends to
zero linearly. Finally, for sufficiently small coupling, the sum of
the spectrum with itself is an interval. This last result provides a
rigorous explanation of a phenomenon for the Fibonacci square
lattice discovered numerically by Even-Dar Mandel and Lifshitz.
\end{abstract}

\maketitle

\section{Introduction}

The Fibonacci Hamiltonian is a central model in the study of
electronic properties of one-dimensional quasicrystals. It is
given by the following bounded self-adjoint operator in
$\ell^2(\Z)$,
$$
[H_{V,\omega} \psi] (n) = \psi(n+1) + \psi(n-1) + V
\chi_{[1-\alpha , 1)}(n\alpha + \omega \!\!\! \mod 1) \psi(n),
$$
where $V > 0$, $\alpha = \frac{\sqrt{5}-1}{2}$, and $\omega \in \T
= \R / \Z$.

This operator family has been studied in many papers since the
early 1980's and numerous fundamental results are known. Let us
recall some of them and refer the reader to the survey articles
\cite{D00,D07,S95} for additional information.

The spectrum is easily seen to be independent of $\omega$ and may
therefore be denoted by $\Sigma_V$. That is, $\sigma(H_{V,\omega})
= \Sigma_V$ for every $\omega \in \T$. Indeed, this follows
quickly from the minimality of the irrational rotation by $\alpha$
and strong operator convergence. It was shown by S\"ut\H{o} that
$\Sigma_V$ has zero Lebesgue measure for every $V > 0$; see
\cite{S89}. Moreover, it is compact (since it is the spectrum of a
bounded operator) and perfect (because the irrational rotation by
$\alpha$ is ergodic). Thus, $\Sigma_V$ is a zero-measure Cantor
set. This result was recently strengthened by Cantat \cite{Can}
who showed that the Hausdorff dimension of $\Sigma_V$ lies
strictly between zero and one.

Naturally, one is interested in fractal properties of $\Sigma_V$,
such as its dimension, thickness, and denseness.  While such a study
is well-motivated from a purely mathematical perspective, we want to
point out that there is significant additional interest in these
quantities. In particular, it has recently been realized that the
fractal dimension of the spectrum is intimately related with the
long-time asymptotics of the solution to the associated
time-dependent Schr\"odinger equation, that is, $i \partial_t \phi =
H_{V,\omega} \phi$; see \cite{DEGT}.

Fractal properties of $\Sigma_V$ are by now well understood for large values of
$V$. Work of Casdagli \cite{Cas} and S\"ut\H{o} \cite{S87} shows
that for $V \ge 16$, $\Sigma_V$ is a dynamically defined Cantor
set. It follows from this result that the Hausdorff dimension and
the upper and lower box counting dimension of $\Sigma_V$ all
coincide; let us denote this common value by $\dim
\Sigma_V$. Using this result, Damanik, Embree, Gorodetski, and
Tcheremchantsev have shown upper and lower bounds for the
dimension; see \cite{DEGT}. A particular consequence of these
bounds is the identification of the asymptotic behavior of the
dimension as $V$ tends to infinity:
$$
\lim_{V \to \infty} \dim \Sigma_V \cdot \log V = \log (1 + \sqrt{2}).
$$
The paper \cite{DEGT} also discusses some of the implications for
the dynamics of the Schr\"odinger equation; let us mention
\cite{DT07,DT08} for further recent advances in this direction for
the strongly coupled Fibonacci Hamiltonian.

By contrast, hardly anything about $\Sigma_V$ (beyond it having
 Hausdorff dimension strictly between zero and one) is known for small values of $V$.
The largeness of $V$ enters the proofs of the existing results in
critical ways. Consequently, these proofs indeed break down once
the largeness assumption is dropped. The purpose of this note is
to announce results about $\Sigma_V$ for $V$ sufficiently small
that are shown by completely different methods. We will indicate
some of the main ideas that are used to prove these results, but
we defer full details to a future publication.

We would like to emphasize that quantitative properties of regular
Cantor sets such as thickness and denseness are widely used in
dynamical systems (see \cite{N70, N79}, \cite{PT}, \cite{M}), found
an application in number theory (see \cite{H}), but to the best of
our knowledge, these kinds of techniques have never been used before
in the context of mathematical physics.

\section{Statement of the Results}

In this section we describe our results for small coupling $V$.
Clearly, as $V$ approaches zero, $H_{V,\omega}$ approaches the
free Schr\"odinger operator
$$
[H_0 \psi] (n) = \psi(n+1) + \psi(n-1),
$$
which is a well-studied object whose spectral properties are
completely understood. In particular, the spectrum of $H_0$ is
given by the interval $[-2,2]$. It is natural to ask which
spectral features of $H_{V,\omega}$ approach those of $H_0$. It
follows from S\"ut\H{o}'s result that the Lebesgue measure of the
spectrum does not extend continuously to the case $V=0$. Given
this situation, one would at least hope that the dimension of the
spectrum is continuous at $V = 0$.

It was shown by us in \cite{DG} (and independently by Cantat \cite{Can}) that $\Sigma_V$ is a dynamically
defined Cantor set for $V>0$ sufficiently small (i.e., the small
coupling counterpart to Casdagli's result at large coupling). A
consequence of this is the equality of Hausdorff dimension and
upper and lower box counting dimensions of $\Sigma_V$ in this
coupling constant regime.  Our first result shows that the dimension of the
spectrum indeed extends continuously to $V=0$.

\begin{Thm}\label{t.1}
We have
$$
\lim_{V \to 0} \dim \Sigma_V = 1.
$$
More precisely, there are constants $C_1, C_2 > 0$ such that
$$
1 - C_1 V \le \dim \Sigma_V \le 1 - C_2 V
$$
for $V > 0$ sufficiently small.
\end{Thm}

Theorem~\ref{t.1} is a consequence of a connection between the
Hausdorff dimension of a Cantor set and its denseness and
thickness, along with estimates for the latter quantities. Since
these notions and connections may be less familiar to at least a
part of our intended audience, let us recall the definitions and
some of the main results; an excellent general reference in this
context is \cite{PT}.

Let $C \subset \mathbb{R}$ be a Cantor set and denote by $I$ its
convex hull. Any connected component of $I\backslash C$ is called a
\emph{gap} of $C$. A \emph{presentation} of $C$ is given by an
ordering $\mathcal{U} = \{U_n\}_{n \ge 1}$ of the gaps of $C$. If $u
\in C$ is a boundary point of a gap $U$ of $C$, we denote by $K$ the
connected component of $I\backslash (U_1\cup U_2\cup\ldots \cup
U_n)$ (with $n$ chosen so that $U_n = U$) that contains $u$ and
write
$$
\tau(C, \mathcal{U}, u)=\frac{|K|}{|U|}.
$$

With this notation, the \emph{thickness} $\tau(C)$ and the
\emph{denseness} $\theta(C)$ of $C$ are given by
$$
\tau(C) = \sup_{\mathcal{U}} \inf_{u} \tau(C, \mathcal{U}, u),
\qquad \theta(C) = \inf_{\mathcal{U}} \sup_{u} \tau(C,
\mathcal{U}, u),
$$
and they are related to the Hausdorff dimension of $C$ by the
following inequalities (cf.~\cite[Section~4.2]{PT}),
$$
\dim_\mathrm{H} C \ge \frac{\log
2}{\log(2+ \frac{1}{\tau(C)})}, \qquad \dim_\mathrm{H} C \le
\frac{\log 2}{\log(2 + \frac{1}{\theta(C)})}.
$$

Due to these inequalities, Theorem \ref{t.1} is a consequence of the following result:

\begin{Thm}\label{t.2}
We have
$$
\lim_{V \to 0}\tau(\Sigma_V) = \infty.
$$
More precisely, there are constants $C_3, C_4 > 0$ such that
$$
C_3 V^{-1} \le \tau(\Sigma_V)\le \theta(\Sigma_V) \le  C_4 V^{-1}
$$
for $V > 0$ sufficiently small.
\end{Thm}

Bovier and Ghez described in their 1995 paper \cite{BG} the
then-state of the art concerning mathematically rigorous results
for Schr\"odinger operators in $\ell^2(\Z)$ with potentials
generated by primitive substitutions. The Fibonacci Hamiltonian
belongs to this class; more precisely, it is in many ways the most
important example within this class of models. One of the more
spectacular discoveries is that, in this class of models, the
spectrum jumps from being an interval for coupling $V = 0$ to
being a zero-measure Cantor set for coupling $V > 0$. That is, as
the potential is turned on, a dense set of gaps opens immediately
(and the complement of these gaps has zero Lebesgue measure). It
is natural to ask about the size of these gaps, which can in fact
be parametrized by a canonical countable set of gap labels; see
\cite{BBG92}. For some examples, these gap openings were studied
in \cite{B} and \cite{BBG91}. However, for the important Fibonacci
case, the problem remained open. In fact, Bovier and Ghez write on
p.~2321 of \cite{BG}: \textit{It is a quite perplexing feature
that even in the simplest case of all, the golden Fibonacci
sequence, the opening of the gaps at small coupling is not known!}

Our next result resolves this issue completely and shows that, in
the Fibonacci case, all gaps open linearly:

\begin{Thm}\label{t.3}
For $V > 0$ sufficiently small, the boundary points of a gap in
the spectrum $\Sigma_V$ depend smoothly on the coupling constant
$V$. Moreover, given any such one-parameter family $\{U_V\}_{V >
0}$, where $U_V$ is a gap of $\Sigma_V$ and the boundary points of
$U_V$ depend smoothly on $V$, we have that
$$
\lim_{V\to 0} \frac{|U_V|}{V}
$$
exists and belongs to $(0,\infty)$.
\end{Thm}

Our next result concerns the sum set $\Sigma_V + \Sigma_V = \{ E_1
+ E_2 : E_1 , E_2 \in \Sigma_V\}$. This set is equal to the
spectrum of the so-called square Fibonacci Hamiltonian. Here, one
considers the Schr\"odinger operator
\begin{align*}
[H^{(2)}_V \psi] (m,n) = & \psi(m+1,n) + \psi(m-1,n) + \psi(m,n+1)
+ \psi(m,n-1) + \\
& + V \left( \chi_{[1-\alpha , 1)}(m\alpha \!\!\! \mod 1) +
\chi_{[1-\alpha , 1)}(n\alpha \!\!\! \mod 1) \right) \psi(m,n)
\end{align*}
in $\ell^2(\Z^2)$. The theory of tensor products of Hilbert spaces
and operators then implies that $\sigma(H^{(2)}_V) = \Sigma_V +
\Sigma_V$. This operator and its spectrum have been studied
numerically and heuristically by Even-Dar Mandel and Lifshitz in a
series of papers \cite{EL06, EL07, EL08}. Their study suggested
that at small coupling, the spectrum of $\Sigma_V + \Sigma_V$ is
\text{not} a Cantor set; quite on the contrary, it has no gaps at
all.

Our final theorem confirms this observation:

\begin{Thm}\label{t.4}
For $V > 0$ sufficiently small, we have that $\Sigma_V + \Sigma_V$
is an interval.
\end{Thm}

Notice that Theorem \ref{t.4} is a direct consequence of Theorem
\ref{t.2} and the famous Gap Lemma, which was used by Newhouse to
construct persistent tangencies and generic diffeomorphisms with
infinite number of attractors (the so-called ``Newhouse
phenomenon"):

\begin{GL}[Newhouse \cite{N70}]\label{t.gaplemma}
If $C_1, C_2\subset \mathbb{R}^1$ are Cantor sets such that
$$
\tau(C_1) \cdot \tau(C_2)>1,
$$
then either one of these sets is contained entirely in a gap of the
other set, or $C_1\cap C_2\ne \emptyset$.
\end{GL}

\section{Comments on the Proofs}

We will exploit the following relation between the spectrum of the
Fibonacci Hamiltonian and the dynamical system known as the trace
map. By the trace map we mean the following polynomial map:
$$
T:\mathbb{R}^3\to \mathbb{R}^3, \ \ \ T(x, y, z)= (2xy-z, x, y).
$$
The trace map preserves the following quantity (called Fricke-Vogt invariant):
$$
I(x,y,z) = x^2 + y^2 + z^2 - 2xyz - 1.
$$
We denote by $S_V$ the surface $\left\{I(x, y,
z)=\frac{V^2}{4}\right\}$ invariant under $T$.

The relation between the trace map and the spectrum of the
Fibonacci Hamiltonian is given by the following statement.

\bthm[S\"ut\H{o} \cite{S87}]\label{spectrum}  An energy $E \in \R$
belongs to the spectrum $\Sigma_V$ of $H_{V,\omega}$ if and only
if the positive semiorbit of the point $(\frac{E-V}{2},
\frac{E}{2}, 1)$ under iterates of the trace map $T$ is bounded.
\ethm

Denote by $\ell_V$ the line
$$
\ell_V = \left\{ \left(\frac{E-V}{2}, \frac{E}{2}, 1 \right) : E
\in \Bbb{R} \right\}.
$$
It is easy to check that $\ell_V \subset S_V$. It is known that
the set
$$
\Lambda_V = \left\{p\in S_V : \ \mathcal{O}_T(p) \ \text{\rm is \
bounded} \right\}
$$
is a hyperbolic set for all $V\ne 0$, see \cite{Cas}, \cite{DG},
\cite{Can}. Moreover, its stable manifolds are transversal to
$\ell_V$ for small $V$. Therefore the set $W^s(\Lambda_V)\cap
\ell_V$ is affinely equivalent to $\Sigma_V$. In particular,
$\Sigma_V$ is a dynamically defined Cantor set for small values of
$V$.

The surface $S_0$ is the so-called Cayley cubic, it has four conic
singularities and can be represented as a union of a two dimensional
sphere (with four conic singularities) and four unbounded
components. The restriction of the trace map to the sphere is a
pseudo-Anosov map (a factor of a hyperbolic map of a two-torus), and
its Markov partition can be presented explicitly (see \cite{Cas} or
\cite{DG}). For small values of $V$, the map $T : S_V\to S_V$
``inherits'' the hyperbolicity of this pseudo-Anosov map everywhere
away from singularities. The dynamics near the singularities needs
to be considered separately. Due to the symmetries of the trace map,
it is enough to consider the dynamics of $T$ near one of the
singularities, say, near the point $p=(1, 1, 1)$. The set $Per_2(T)$
of periodic orbits of period two is a smooth curve that contains the
point $p$ and intersects $S_V$ at two points (denote them by $p_1$
and $p_2$) for $V > 0$. Also, this curve is a normally hyperbolic
manifold, and we can use the normally hyperbolic theory (see
\cite{HPS}, \cite{PSW}) to study the behavior of $T$ in a small
neighborhood of $p$; see \cite{DG} for details. In particular, since
the distance between $p_1$ and $p_2$ is of order $V$, and the gaps
in the spectrum are formed by the points of intersection of
$W^{ss}(p_1)$ and $W^{ss}(p_2)$ with the line $\ell_V$, the size of
a given gap is of order $V$ as $V\to 0$, which implies Theorem
\ref{t.3}.

In order to estimate the thickness (and the denseness) of the
spectrum $\Sigma_V$, we notice first that the Markov partition for
$T:S_0\to S_0$ can be continuously extended to a Markov partition
for $T:\Lambda_V\to \Lambda_V$, and while the size of the elements
of these Markov partitions remains bounded, the size of the
distance between them is of order $V$. The natural approach now is
to use the distortion property (see, e.g., \cite{PT}) to show that
for the iterated Markov partition, the ratio of the distance
between the elements to the size of an element is of the same
order. The main technical problem here is again the dynamics of
the trace map near the singularities, since the curvature of $S_V$
is very large there for small $V$. Nevertheless, one can still
estimate the distortion that is obtained during a transition
through a neighborhood of a singularity and prove boundedness of
the distortion for arbitrarily large iterates of the trace map.
This implies Theorem~\ref{t.2}.

\end{document}